\documentclass{article}

\usepackage{amssymb, amsmath, graphicx}
\usepackage{fullpage}

\begin{document}

\title{A Linear-algebraic Proof of Hilbert's Ternary Quartic Theorem}\markright{Hilbert's Ternary Quartic Theorem}
\author{Anatolii Grinshpan and Hugo J.~Woerdeman}


\date{}
\newtheorem{theorem}{Theorem}
\newtheorem{defn}[theorem]{Definition}
\newtheorem{corollary}[theorem]{Corollary}
\newtheorem{lemma}[theorem]{Lemma}
\newtheorem{remark}[theorem]{Remark}
\newtheorem{prob}[theorem]{Problem}
\newtheorem{example}[theorem]{Example}
\newtheorem{proposition}[theorem]{Proposition}
\numberwithin{equation}{section}

\maketitle

\begin{abstract}
Hilbert's ternary quartic theorem states that every nonnegative degree 4 homogeneous polynomial in three variables 
can be written as a sum of three squares of homogeneous quadratic polynomials.
We give a linear-algebraic approach to Hilbert's theorem by showing that a structured cone of positive semidefinite matrices is generated by rank 1 elements.
\end{abstract}

\section{Introduction.}\label{sec:intro}

A homogeneous polynomial $$p(x,y,z)=\sum_{i+j+k=4} p_{ijk}x^iy^jz^k$$ in three variables of degree 4 is called a ternary quartic. Hilbert's classical theorem \cite{Hilbert}, dating back to 1888, states that every ternary quartic that takes only nonnegative values, i.e., such that $p(x,y,z) \ge 0$ for all $x,y,z \in{\mathbb R},$ can be written as a sum of three squares of homogeneous quadratic polynomials. This theorem stood as a precursor of Hilbert's 17th problem and subsequent development, and to this day attracts a lot of attention.
Detailed expositions can be found in \cite{Rudin} and \cite{Swan}.

One distinguishes two parts to Hilbert's theorem: the existence of a representation
as a sum of squares (qualitative part) and the assertion that at most three squares suffice (quantitative part). Hilbert's original proof, cast in modern form, takes roots in advanced topology and algebraic geometry. 
Many attempts have been made in search of more elementary proofs.   
In 1977, Choi and Lam \cite{ChoiLam} gave an elementary proof of the qualitative part, based on properties of extremal positive semidefinite forms. In 2004, Pfister \cite{Pfister} gave a different elementary proof, which was constructive. 
New approaches to Hilbert's theorem were developed in \cite{PS} and \cite{PRS}. 
But no simple elementary explanation of the quantitative part has been found.
Very recently, Hilbert's theorem has been considered from a new general perspective, in the
framework of nonnegative quadratic forms on projective real varieties \cite{BPSV, BSV}. 

In this note we would like to offer a new elementary proof of the qualitative part of Hilbert's theorem.
Our approach uses linear algebra and convex geometry.

\section{The PSD$_6$ cone and Hilbert's theorem.}\label{sec:cone}

We begin with a few preliminary facts. As general sources, we refer the reader to \cite{HornJ}
for background on matrix theory and to \cite{Barvinok, B, BV, Webster} for background on convex geometry.

Let $S_n$ be the vector space of $n\times n$ real symmetric matrices $A=A^\top$ (the superscript $\top$ denotes the transpose). The dimension of $S_n$ is
$n(n+1)/2$. The scalar product of two symmetric matrices $A=(a_{ij}), B=(b_{ij})$ in $S_n$ is defined by
$$\langle A, B\rangle=\sum_{i,j=1}^na_{ij}b_{ij}={\rm tr}(AB),$$
where tr, the trace, is the sum of diagonal elements of a matrix. Equipped with the scalar product, $S_n$ becomes
a Euclidean space.
Every hyperplane in $S_n$ is of the form $$\{ X \in S_n:\ {\rm tr} (XC)=h \},$$ where $C \in S_n$ is a nonzero matrix and $h \in {\mathbb R}$. Two subsets of $S_n$ are said to be {\it strictly separated\rm}  by the hyperplane $ {\rm tr} (XC) =h$ if one is contained in the open half-space $\{ X :  {\rm tr} (XC) <h \}$ and the other in the open half-space $\{ X : {\rm tr} (XC) > h \}$. Two disjoint closed convex sets in a Euclidean space can be strictly separated by a hyperplane if their vector difference is closed \cite[Proposition 1.5.3]{B}.

%

A matrix $A\in S_n$ is said to be positive semidefinite if $\langle Ax,\ x\rangle\ge0$ for all vectors $x\in\mathbb R^n$. Equivalently, $A$ is positive semidefinite if there is a matrix $B$ such that $A=B B^\top$. 
In particular, if $A$ is of rank $k$, then $B$ can be chosen of size $n\times k$.

The set PSD$_n$ of all $n\times n$ positive semidefinite matrices is closed under addition and nonnegative scaling. 
Such a set is said to be a cone or, more precisely, a convex cone. The cone PSD$_n$ is pointed (i.e., contains no lines) and closed.

A ray of a cone generated by its (nonzero) element consists of all nonnegative multiples of the element. 
A ray of a cone is called extreme if it cannot be expressed as a nonnegative linear combination of other rays. Minkowski's theorem for cones asserts that every ray of a closed pointed (convex) cone is a nonnegative linear combination of its extreme rays \cite[Sections II.3 and II.8]{Barvinok}. In particular, Minkowski's theorem applies to PSD$_n$.

For every $X\in S_n$, the condition ${\rm tr}(XY)\ge0$, for all $Y\in$\,PSD$_n$, is equivalent to $X\in$\,PSD$_n$. This is known as self-duality of the PSD$_n$ cone \cite[Section 2.6.1]{BV}.

The connection to Hilbert's theorem can now be explained.
A polynomial $p(x,y,z)$ is a sum of squares of homogeneous quadratic polynomials if and only if it can be represented in the form
\begin{equation}\label{pxyz} p(x,y,z) =  \begin{bmatrix} x^2 & xy & xz & y^2 & yz & z^2 \end{bmatrix} A \begin{bmatrix} x^2 \cr xy \cr xz \cr y^2 \cr yz \cr z^2 \end{bmatrix}, \end{equation}
where $A \in {\rm PSD}_6$. Indeed, if $A=\sum_{i=1}^k a_i a_i^\top$, where $a_i\in {\mathbb R}^6$, then \eqref{pxyz} turns into a desired sum-of-squares representation:
$$ p(x,y,z) = \sum_{i=1}^k (v^\top a_i)^2,\qquad v^\top= \begin{bmatrix} x^2 & xy & xz & y^2 & yz & z^2 \end{bmatrix}.$$
In fact, representation \eqref{pxyz} is easy to obtain if we merely require $A$ to be symmetric and thus drop the positive semidefiniteness condition. One such choice is given by
$$ A_0 = \frac12\begin{bmatrix} 2p_{400} &  p_{310} & p_{301} & 0 & p_{211} & 0\cr
 p_{310} & 2p_{220} & 0 & p_{130} & p_{121} & 0 \cr
 p_{301} & 0 &  2p_{202} & 0 & p_{112} & p_{103}  \cr 
 0 & p_{130} & 0 &  2p_{040} & p_{031} & 0 \cr
 p_{211} & p_{121} &  p_{112} & p_{ 031} & 2p_{022} & p_{013} \cr
 0 & 0 & p_{103} & 0 & p_{013} & 2p_{004} 
 \end{bmatrix}. $$
Moreover, if ${\mathcal W}$ is the subspace of $S_6$ consisting of the matrices 
$$\begin{bmatrix} 
0 & 0 & 0 & w_1 & w_2 & w_3 \cr 
0 & -2w_1 & -w_2 & 0 & w_4 & w_5 \cr
0 & -w_2 & -2w_3 & -w_4 & -w_5 & -0 \cr 
w_1 & 0 & -w_4 & 0 & 0 & w_6 \cr 
w_2 & w_4 & -w_5 & 0 & -2w_6 & 0\cr 
w_3 & w_5 & 0 & w_6 & 0 & 0
\end{bmatrix},\qquad w_1, \ldots, w_6 \in {\mathbb R}, $$
then \eqref{pxyz} holds if and only if $A\in A_0 + {\mathcal W}$. Thus it suffices to show that 
\begin{equation}\label{ne} (A_0 + {\mathcal W}) \cap {\rm PSD}_6 \neq \varnothing . \end{equation}
It turns out that condition \eqref{ne} holds if and only if there is no hyperplane strictly separating the convex sets $A_0+{\mathcal W}$ and ${\rm PSD}_6$. The two possible scenarios are illustrated in Figure 1.
\begin{figure}[h]
\begin{center}
\includegraphics[scale=.58]{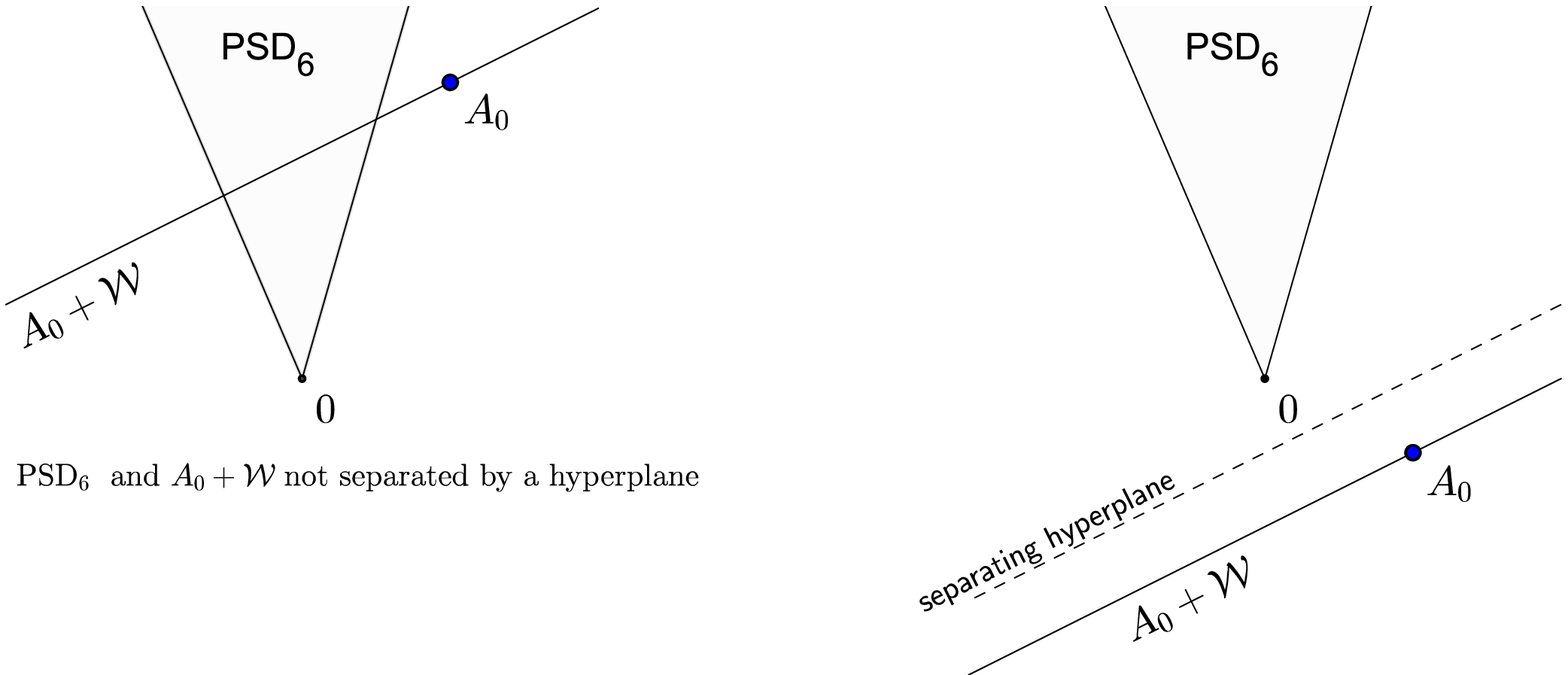}
\end{center}
\caption{$A_0+\mathcal W$ and PSD$_6$ cannot be strictly separated if and only if $(A_0+\mathcal W)\cap{\rm PSD}_6\ne\varnothing$.}
\end{figure}

A hyperplane $\{ X \in S_6:\ {\rm tr} (XC) =h \}$ can be disjoint from $A_0 + {\mathcal W}$ only if $C$ belongs to ${\mathcal W}^\perp$, the orthogonal complement of $\mathcal W$ in $S_6$. 
Note that the subspace ${\mathcal W}^\perp$ consists of all real symmetric matrices $(a_{ij})_{i,j=1}^6$ such that
\begin{equation}\label{C}
\begin{aligned}
&a_{14}=a_{22}, && a_{15}=a_{23}, && a_{16}=a_{33},\\
&a_{25}=a_{34}, && a_{26}=a_{35}, && a_{46} = a_{55}.
\end{aligned}
\end{equation}
Since ${\mathcal W}^\perp$ contains a positive definite matrix, the intersection of ${\rm PSD}_6$ and ${\mathcal W}$ contains only the zero matrix. By \cite[Proposition 1.4.14]{B} the vector difference ${\rm PSD}_6 - {\mathcal W}$ is closed and consequently there is a hyperplane strictly separating ${\rm PSD}_6$ and $A_0+{\mathcal W}$ \cite[Proposition 1.5.3]{B}.
Now if PSD$_6$ is contained in the closed half-space $$\{ X \in S_6:\ {\rm tr} (XC) \ge h \},$$ then
$h\le 0$ and $C \in {\rm PSD}_6$, by the self-duality of ${\rm PSD}_6$.
Therefore \eqref{ne} holds if and only if 
\begin{equation}\label{ne2} C \in {\rm PSD}_6 \cap {\mathcal W}^\perp\quad {\rm implies}\quad  {\rm tr} (A_0C) \ge 0 . \end{equation}
The key in proving implication \eqref{ne2} is our main result, which we now state.

\begin{theorem}\label{main}
Let $\ \mathcal C$ be the cone of positive semidefinite matrices in $S_6$ satisfying \eqref{C}.
Then every extreme ray of $\ \mathcal C$ is generated by a rank 1 matrix $vv^\top$, where
\begin{equation}\label{v} v^\top= v(x,y,z)^\top:= \begin{bmatrix} x^2 & xy & xz & y^2 & yz & z^2 \end{bmatrix}, \end{equation} for some $ x,\ y,\ z \in\mathbb R$. 
Thus every element of $\ {\mathcal C}$ is a nonnegative linear combination of matrices $vv^\top$. 
\end{theorem}

The second assertion of Theorem \ref{main} follows by Minkowski's theorem.
Note that if $p(x,y,z) $ takes only nonnegative values, we obtain that 
$${\rm tr} (A_0 vv^\top ) ={\rm tr} (v^\top A_0 v ) = p(x,y,z) \ge 0.$$ 
By Theorem \ref{main}, each element $C$ of the cone ${\mathcal C} =  {\rm PSD}_6 \cap {\mathcal W}^\perp$ is of the form $C = \sum_i \rho_i v_i v_i^\top$ with $\rho_i \ge 0$ and $v_i=v(x_i,y_i,z_i)$, and thus 
$${\rm tr} (A_0 C ) = \sum_i \rho_i p(x_i , y_i , z_i ) \ge 0.$$ Consequently, Theorem \ref{main} proves that a ternary quartic that takes only nonnegative values is a sum of squares.

\section{Proof of the Theorem.}\label{sec:proof}

We now prove Theorem \ref{main}. The argument hinges on the following lemma.
\begin{lemma} \label{vectors}
Let 
$u,v,w,y \in {\mathbb R}^2$ be such that $\langle u,v \rangle = \langle w,y \rangle$. Then there exists a rotation $R$ with the property that 
$$R\begin{bmatrix} u & v & w & y \end{bmatrix}=\begin{bmatrix} u_1 & v_1 & w_1 & y_1 \cr u_2 & v_2 & w_2 & y_2 \end{bmatrix}$$
satisfies $u_1v_1 = w_1 y_1$ and $u_2 v_2 = w_2 y_2 $. 
\end{lemma}
We alert the reader that the subscripts in Lemma \ref{vectors} are used to indicate the components of the rotated, not original, vectors. The same convention applies further below in the proof of Theorem \ref{main}.

{\it Proof.\rm} Though the statement is about vectors in $u,v,w,y \in{\mathbb R}^2$, it is convenient to treat them as complex numbers $\alpha , \beta , \gamma , \delta$,
$$ \begin{bmatrix} \alpha & \beta & \gamma & \delta \end{bmatrix} =  \begin{bmatrix} 1 & i \end{bmatrix}  \begin{bmatrix} u & v & w & y \end{bmatrix} . $$ 
The task  amounts to choosing an angle of rotation $\theta$ so that the rotated complex numbers satisfy
\begin{equation}\label{rot}
\begin{aligned}
&{\rm Re} (e^{i\theta}\alpha) {\rm Re}(e^{i\theta} \beta) - {\rm Re}(e^{i\theta}\gamma)  {\rm Re}(e^{i\theta} \delta)=u_1v_1-w_1y_1=0\\ 
&{\rm Im}(e^{i\theta}\alpha) {\rm Im}(e^{i\theta} \beta)- {\rm Im}(e^{i\theta}\gamma)  {\rm Im}(e^{i\theta}  \delta)=u_2v_2-w_2y_2=0.
\end{aligned}
\end{equation}
But since the scalar products $\langle u,v \rangle,\ \langle w,y \rangle$ are rotation invariant, i.e.,
\begin{align*}
&{\rm Re} (e^{i\theta}\alpha) {\rm Re}(e^{i\theta} \beta)+{\rm Im} (e^{i\theta}\alpha) {\rm Im}(e^{i\theta} \beta)={\rm Re}(e^{i\theta}\alpha\,\overline{e^{i\theta}\beta})={\rm Re}(\alpha\,\overline{\beta})\\
&{\rm Re} (e^{i\theta}\gamma) {\rm Re}(e^{i\theta} \delta)+{\rm Im} (e^{i\theta}\gamma) {\rm Im}(e^{i\theta} \delta)={\rm Re}(e^{i\theta}\gamma\,\overline{e^{i\theta}\delta})={\rm Re}(\gamma\,\overline{\delta}),
\end{align*}
the assumption $\langle u,v \rangle=\langle w,y \rangle$ means that the sum of numbers
$${\rm Re} (e^{i\theta}\alpha) {\rm Re}(e^{i\theta} \beta)-{\rm Re} (e^{i\theta}\gamma) {\rm Re}(e^{i\theta} \delta),\quad 
{\rm Im}(e^{i\theta}\alpha) {\rm Im}(e^{i\theta} \beta)- {\rm Im}(e^{i\theta}\gamma)  {\rm Im}(e^{i\theta}  \delta)$$
is zero. Thus, to obtain \eqref{rot} it suffices to choose $\ \theta\ $ so that the difference of these numbers is also zero, which reduces to 
$${\rm Re}(e^{2i\theta} (\alpha \beta - \gamma \delta))=0.$$
The latter is clearly possible and this establishes Lemma \ref{vectors}.
\hfill $\Box$

{\it Proof of Theorem \ref{main}.\rm} Let us start by observing that a rank 1 element in the cone ${\mathcal C}$ is necessarily of the form $vv^\top$, with $v$ as in \eqref{v}. Indeed, if $A=(a_{ij})_{i,j=1}^6 \in {\mathcal C}$ is of rank 1, then its diagonal entries are nonnegative, and we can introduce 
$$ x=\sqrt[4]{a_{11}},\quad  y =   {\rm sign}(a_{12})\sqrt[4]{a_{44}},\quad z = {\rm sign}(a_{13}) \sqrt[4]{a_{66}}.$$ 
As $A$ is of rank 1 (and symmetric), we obtain that $a_{14}^2 = a_{11}a_{44} = x^4y^4$. Since $a_{14}=a_{22}\ge0$, the equalities $a_{22}=a_{14}=x^2y^2$ follow. Similarly, $a_{33}=a_{16} = x^2z^2$ and $a_{55}= a_{46}=y^2z^2$.  Next, $a_{12}^2 = a_{11}a_{22} = x^6 y^2$, and since $x\ge 0$ and ${\rm sign} (a_{12}) = {\rm sign} (y)$, we have $a_{12}=x^3y$. Continuing in this way, we obtain expressions for all entries of $A$ in terms of $x,y,$ and $z$, and $A=vv^\top$ follows.  

Prior to analyzing the rank 2 elements of ${\mathcal C}$, let us make some useful observations. 
As before, let $S_k$ denote the vector space of $k\times k$ real symmetric matrices equipped with the scalar product $\langle A, B \rangle =  {\rm tr}(AB)$. Writing a rank $k$ element $A$ of $\mathcal C$ as \begin{equation}\label{A} A = \begin{bmatrix} a^\top \cr b^\top 
\cr c^\top \cr d^\top \cr e^\top \cr f^\top \end{bmatrix} \begin{bmatrix} a & b & c & d & e & f \end{bmatrix},\qquad a,\ b,\ c,\ d,\ e,\ f \in\mathbb R^k, \end{equation} we have that the linear constraints \eqref{C} are equivalent to the conditions
\begin{equation}\label{E} E_i \in \left({\rm span}\ I_k \right)^\perp,\qquad i=1,2,3,4,5,6, \end{equation} where
$I_k$ is the $k\times k$ identity matrix and
\begin{align*}
&E_1= ad^\top+ da^\top-2bb^\top,\quad &&E_2 = af^\top+ fa^\top - 2cc^\top,\\
&E_3= bc^\top+ cb^\top-ea^\top-ae^\top,\quad &&E_4= be^\top+eb^\top-cd^\top-dc^\top,\\
&E_5= ce^\top+ec^\top-bf^\top- fb^\top,\quad &&E_6= df^\top+ fd^\top-2ee^\top
\end{align*} 
are matrices in $S_k$. When $k \ge 4$, there exists a nonzero $F\in S_k$ orthogonal to $I_k,\,E_1,\ldots,\,E_6,$ as the dimension of $S_k$ is $k(k+1)/2 \ge 10$.
Consequently, for small $\varepsilon,$ $A$ is the average of distinct points
$$\begin{bmatrix} a^\top \cr b^\top
\cr c^\top \cr d^\top \cr e^\top \cr f^\top \end{bmatrix} (I\pm \varepsilon F) \begin{bmatrix} a & b & c & d& e & f \end{bmatrix}\in\mathcal C$$  
and does not generate an extreme ray.

Now consider the rank 2 case. Let $A$ be as in \eqref{A} with $k=2$. Condition \eqref{E} means that
\begin{align*}
&\langle a,d \rangle = \langle b,b \rangle, &&\langle a,f \rangle = \langle c,c \rangle, &&\langle b,c \rangle = \langle a,e \rangle,\\
&\langle b,e \rangle = \langle c,d \rangle, &&\langle c,e \rangle = \langle b,f \rangle, &&\langle d,f \rangle = \langle e,e \rangle.
\end{align*} 
Let us assume that no vectors among $a,\ b,\ c,\ d,\ e,\ f$ are multiples of each other. As $\langle b,c \rangle = \langle a,e \rangle,$ 
by Lemma \ref{vectors} there exists a rotation $R$ so that 
$$R\begin{bmatrix} a & e & b & c \end{bmatrix}=\begin{bmatrix} a_1 & e_1 & b_1 & c_1 \cr a_2 & e_2 & b_2 & c_2 \end{bmatrix}$$
satisfies $a_1e_1 = b_1 c_1$ and $a_2 e_2 = b_2 c_2 $.
Now write 
$$ A = \begin{bmatrix} a^\top \cr b^\top
\cr c^\top \cr d^\top \cr e^\top \cr f^\top \end{bmatrix} R^\top R \begin{bmatrix} a & b & c & d& e & f \end{bmatrix}
= \begin{bmatrix} a_1 & a_2 \cr b_1 & b_2 
\cr c_1 & c_2 \cr d_1 & d_2 \cr e_1 & e_2 \cr f_1 & f_2 \end{bmatrix} \begin{bmatrix} a_1 & b_1 & c_1 & d_1& e_1 & f_1 \cr a_2 & b_2 & c_2 & d_2& e_2 & f_2\end{bmatrix} .$$
Assuming that $a_1,\ a_2\ne 0,$ let $\tilde d = \begin{bmatrix}\frac{b_1^2}{a_1} & \frac{b_2^2}{a_2}\end{bmatrix}^\top\!\!\!,$
so that $\langle b, e \rangle = \langle c,\tilde d \rangle$ and $\langle a,\tilde d \rangle = \langle b,b \rangle.$
This yields $$ \langle a,d-\tilde d \rangle = 0 = \langle c, d-\tilde d \rangle, $$
and as $a$ and $c$ are linearly independent, we get that $d=\tilde d$, yielding $d_1=\frac{b_1^2}{a_1}$ and $d_2=\frac{b_2^2}{a_2}$. Similarly, we find that $f_1=\frac{c_1^2}{a_1}$ and $f_2=\frac{c_2^2}{a_2}$. So letting
\begin{align*}
&x=\sqrt{|a_1|}, &&y=\frac{b_1}{x}, &&z=\frac{c_1}{x}, \\
&\hat x = \sqrt{|a_2|}, &&\hat y=\frac{b_2}{\hat x}, &&\hat z=\frac{c_2}{\hat x}, 
\end{align*}
one easily checks that $A=vv^\top + \hat v \hat v^\top$, where $v$ and $\hat v$ are as in \eqref{v}.

\begin{remark}\label{lincomb}  
\rm In the above reasoning it would have sufficed to have the following equalities from the start:
\begin{align*}
&\langle a,d \rangle = \langle b,b \rangle, &&\langle a,f \rangle = \langle c, c \rangle,\hspace{.5in} \langle b,c \rangle = \langle a,e \rangle, \\
&\langle b,e \rangle = \langle c,d \rangle, &&\alpha \langle f,b \rangle + \beta \langle f,d \rangle = \alpha \langle e,c \rangle + \beta \langle e,e \rangle,   
\end{align*}
for some scalars $\alpha, \beta$ such that $a$ and $\alpha b+\beta d$ are linearly independent. Indeed, the equality $\langle b,c \rangle = \langle a,e \rangle$ gives $a_1e_1 = b_1 c_1$ and $a_2 e_2 = b_2 c_2 $. The equalities $\langle b,e \rangle = \langle c,\tilde d \rangle, \langle a,\tilde d \rangle = \langle b,b \rangle $ then give $d_1=\frac{b_1^2}{a_1}$ and $d_2=\frac{b_2^2}{a_2}$.
To conclude that $f_1=\frac{c_1^2}{a_1}$ and $f_2=\frac{c_2^2}{a_2}$, we use that both $f$ and $\tilde f := \begin{bmatrix}  \frac{c_1^2}{a_1} & \frac{c_2^2}{a_2} 
\end{bmatrix}^\top$ satisfy the conditions $\langle f,a \rangle = \langle c,c\rangle $ and $\alpha \langle f,b \rangle + \beta \langle f,d \rangle = \alpha \langle e,c \rangle + \beta \langle e,e \rangle$. Thus $ \langle a, f-\tilde f \rangle = 0 =
\langle \alpha b +
\beta d , f - \tilde f \rangle$, yielding $f=\tilde f$ as $a$ and $\alpha b + \beta d$ are linearly independent.
\end{remark} 

We have made some generic assumptions in the rank $2$ case, but these can be lifted. When there is some pairwise linear dependence among $a,b,c,d,e,f$ or when $a_1a_2 = 0$, appropriate modifications of the same argument still apply.

Finally, the rank $3$ case remains. Since $\dim{S_3} = 6$, the matrices $E_1, \ldots, E_6$ lie in the 5-dimensional subspace $\{ I_3 \}^\perp$ and so are linearly dependent. Let us assume that $E_6$ is a linear combination of $E_1, \dots , E_5$. Choose a nonzero $F \in \{ E_1, E_2 , E_3, E_4 , I_3 \}^\perp$ and $\varepsilon \neq 0$ so that $I_3-\varepsilon F$ is positive semidefinite with rank 2. Then
$$ B=  \begin{bmatrix} a^\top \cr b^\top 
\cr c^\top \cr d^\top \cr e^\top \cr f^\top \end{bmatrix} (I - \varepsilon F) \begin{bmatrix} a & b & c & d& e & f \end{bmatrix} \ge 0$$ satisfies four of the linear conditions \eqref{C} and a linear combination of the remaining two (as $aE_5+bE_6 \in {\rm span} \{ E_1,E_2,E_3,E_4 \} $). Due to Remark \ref{lincomb} these equalities suffice to show that $B$ does not generate an extreme ray, and therefore $A$ does not either. This finishes the last outstanding case, thus establishing Theorem \ref{main}.
\hfill $\Box$

\section{The number of squares.}\label{sec:number}

As noted in the introduction, Hilbert's result is actually stronger than what we have shown: the sum-of-squares representation can always be chosen to have at most three squares. The known proofs of this fact, including Hilbert's original proof \cite{Hilbert}, are much less elementary \cite{PS, PRS}, and a linear-algebraic argument, if it exists, is yet to be found. 

From the proof of Choi and Lam \cite{ChoiLam} one extracts additional information that every nonnegative ternary quartic is a sum of five squares. Pfister's proof \cite{Pfister} shows that at most four squares suffice.
We note that the four squares conclusion can be reached in a different way, using the geometry of the PSD$_6$ cone.
Namely, by letting $\mathcal A=A_0+\mathcal W$ and $n=6$ in the following lemma of Barvinok.
\begin{lemma}\cite[Chapter II, Lemma 13.6]{Barvinok} Let $\mathcal A$ be an $n$-dimensional affine subspace of $S_n$. If the intersection of $\mathcal A$ with PSD$_n$ is nonempty and bounded, then $\mathcal A$ contains a positive semidefinite matrix of rank at most $n-2$.
\end{lemma}

We conclude with a few words on how one may numerically find a sum-of-squares representation using semidefinite programming (SDP). General references on SDP and convex optimization are \cite{B, BV}. When we let $V_1, \ldots , V_{15}$ be a basis for ${\mathcal W}^\perp$, finding a sum-of-squares representation comes down to finding $A \in {\rm PSD}_6$ with ${\rm tr} (A V_i ) = {\rm tr} (A_0 V_i)=: b_i$, $i=1,\ldots , 15$, which is exactly a feasibility problem in SDP. Choosing a positive definite $C$, one can perform the SDP
$$ \inf {\rm tr} (CA),\qquad {\rm subject \ to }\quad  A\in {\rm PSD}_6,\quad {\rm tr} (A V_i ) =  b_i,\quad i=1,\ldots , 15. $$
In \cite[Section 6]{PSV} it is observed that for random $C$ there is a positive probability to find a rank 3 optimal $A$. Thus a repeated performance of the above SDP with random $C$, ultimately yields a representation as a sum of three squares. Here we accept a solution as having rank at most 3 when its fourth singular value is sufficiently small. It is our experience that this happens after just a few tries. 

{\it Acknowledgment.}\rm
The authors wish to thank the anonymous referees for their thoughtful comments which led to improvement in the presentation. They also thank Benjamin Grossmann
for reading the article and providing helpful feedback. HJW  is supported by Simons Foundation grant 355645.

%

\bigskip

\noindent Department of Mathematics\\
Drexel University \\
Philadelphia, PA 19104\\
\{tolya, hugo\}@math.drexel.edu

\end{document}